%

\input btxmac.tex

\magnification=\magstep1 
\hsize=15truecm
\vsize=22.5truecm
\hoffset=.46truecm

\font\titlefont=cmss12
\font\authorfont=cmr12

\font\amsa=msam10
\def\eop{{\amsa \char"03}}  

\font\tenbb=msbm10     
\font\sevenbb=msbm7
\newfam\bbfam \def\bb{\fam\bbfam\tenbb}  
\textfont\bbfam=\tenbb
\scriptfont\bbfam=\sevenbb

\def\Pt{Pythagorean triple}
\def\intz{{\bb Z}}     
\def\ratq{{\bb Q}}     
\def\x{{\underline x}} 
\def\Int{{\rm Int}}    

\let\congr=\equiv  

\def\[{{\rm[}} \def\]{\/{\rm]}}
\def\({{\rm(}} \def\){\/{\rm)}}

\long\def\profess#1#2\endprofess {\ifdim\lastskip<\medskipamount%
  \removelastskip\penalty-50\medskip\fi
  \noindent{\bf#1\enspace}{\sl#2\par}%
  \ifdim\lastskip<\medskipamount\removelastskip\penalty55\medskip\fi}

\def\proof{\noindent {\it Proof.\enspace}}
\def\endproof{\ \eop\par%
   \ifdim\lastskip<\medskipamount
   \removelastskip\penalty-50\medskip\fi}

\thickmuskip=4.5mu plus 4mu minus 3mu  
\mathsurround=0.5pt

\pretolerance=200  \hyphenpenalty=1000 \exhyphenpenalty=5000

\baselineskip=13pt \parskip=1.5pt \parindent=1.5em
\smallskipamount=4pt plus2pt minus1pt 
\medskipamount=7pt plus3pt minus2pt 
\bigskipamount=15pt plus5pt minus4pt 

\line{\ \hfill \ }
\line{\ {\it To appear in\/} J.~Pure Appl.~Algebra\hfill}
\vskip1truein

\centerline{\titlefont Parametrization of Pythagorean triples}
\smallskip
\centerline{\titlefont by a single triple of polynomials}
\bigskip\medskip

\centerline{\authorfont Sophie Frisch {\tenrm and} Leonid Vaserstein}
\footnote{\vphantom{*}}
{2000 Math.~Subj.~Classification: primary 11D09, secondary 11D85, 11C08, 13F20.}
\bigskip\bigskip\medskip

{\bf Abstract.} 
It is well known that Pythagorean triples can be parametrized by two
triples of polynomials with integer coefficients. We show that no
single triple of polynomials with integer coefficients in any number
of variables is sufficient, but that there exists a parametrization 
of Pythagorean triples by a single triple of integer-valued polynomials.
\bigskip\medskip

The second author has recently studied polynomial parametrizations
of solutions of Diophantine equations \cite{Vas--PPDE}, and the 
first author has wondered if it were possible in some cases to
parametrize by a $k$-tuple of integer-valued polynomials a 
solution set that is not parametrizable by a $k$-tuple of 
polynomials with integer coefficients \cite{Fri--RPP}.
\Pt s provide an example that this is indeed so.

We call a triple of integers $(x,y,z)\in\intz^3$ satisfying
$$ x^2+y^2=z^2$$
a \Pt, and, if $x,y,z>0$, a positive \Pt.

It is well known that every \Pt\ is either of the form 
$$\left(c(a^2-b^2),\; 2cab,\; c(a^2+b^2)\right)$$
or of the form
$$\left(2cab,\; c(a^2-b^2),\; c(a^2+b^2)\right)$$
with $a,b,c\in\intz$, see for instance~\cite{HlaSch79ZT}.

To make precise our usage of the term polynomial parametrization,
consider a set $S\subseteq \intz^k$. Let $f_1,\ldots,f_k$ be a $k$-tuple
of polynomials either in the ring of polynomials with integer 
coefficients in $n$ variables, $\intz[x_1,\ldots,x_n]$, or in
$$\Int(\intz^n)=\{f\in\ratq[x_1,\ldots,x_n]\mid
\forall a\in \intz^n\; f(a)\in\intz\},$$
the ring of integer-valued polynomials in $n$ variables,
for some $n$. In either case $F=(f_1,\ldots,f_k)$
defines a function $F\colon \intz^n\rightarrow \intz^k$.
If $S$ is the image of this function, $S=F(\intz^n)$,
we say that $(f_1,\ldots,f_k)$ parametrizes $S$. We call this a 
parametrization of $S$ by a single $k$-tuple of polynomials.

If $S\subseteq \intz^k$ is the union of the images of finitely many
$k$-tuples of polynomials $F_i=(f_{i1},\ldots,f_{ik})$, 
$S=\bigcup_{i=1}^m F_i(\intz^n)$, we call this a parametrization
of $S$ by a finite number of $k$-tuples of polynomials, and we
distinguish between parametrizations by polynomials with integer
coefficients, meaning $F_i\in (\intz[x_1,\ldots,x_n])^k$ for
all $i$, and by integer-valued polynomials, meaning 
$F_i\in (\Int(\intz^n))^k$.
Unless explicitly specified otherwise, we are using integer parameters,
that is, we let all variables of all polynomials range through the
integers.

We will give a parametrization of the set of \Pt s by a single
triple of integer-valued polynomials in $4$ variables, that is, by 
$(f_1,f_2,f_3)$ with $f_i\in \ratq[x,y,z,w]$ such that
$f_i(x,y,z,w)\in\intz$ whenever $x,y,z,w\in\intz$.
First we will show that it is not possible to parametrize \Pt s 
by a single triple of polynomials with integer coefficients in
any number of variables.

Note, however, that every set of $k$-tuples of integers that is
parametrizable by a single $k$-tuple of integer-valued polynomials 
is parametrizable by a finite number of $k$-tuples of polynomials
with integer coefficients \cite{Fri--RPP}.

\profess{Remark.} 
There do not exist $f,g,h\in \intz[x_1,\ldots,x_n]$ for any $n$
such that $(f,g,h)$ parametrizes the set of \Pt s.
\endprofess

\proof
Suppose $(f,g,h)$ parametrizes the \Pt s. As $\intz[\x]$ is a unique
factorization domain, there exists $d=\gcd(g,h)$ (unique up to sign),
which also divides $f$, since $f^2+g^2=h^2$. 
Let $\varphi=f/d$, $\psi=g/d$ and $\theta=h/d$. 

Then $$\varphi^2= \theta^2 - \psi^2=(\theta+\psi)(\theta-\psi)$$
and $\gcd((\theta+\psi), (\theta-\psi))$ is either $1$ or $2$, but
it can't be $2$, because there exist \Pt s with odd first coordinate
such as $(3,4,5)$. Since $(\theta+\psi)$ and $(\theta-\psi)$ are
co-prime and their product is a square, $(\theta+\psi)$ and
$(\theta-\psi)$ are either both squares, or both $(-1)$ times a square,
and we can get rid of the latter alternative by retro-actively changing
the sign of the polynomial $d$, if necessary.

So there exist polynomials $s$ and $t$ with $(\theta+\psi)=s^2$
and $(\theta-\psi)=t^2$, and therefore 
$$\theta= {s^2+t^2\over 2}\quad \hbox{\rm and}\quad
\psi = {s^2-t^2\over 2}.$$

Since $s^2-t^2 =(s+t)(s-t)$ is divisible by $2$, it is actually 
divisible by $4$, so that $\psi $ is divisible by $2$, which
contradicts the existence of \Pt s with odd second coordinate
such as $(4,3,5)$.
\endproof

%
In a way, it was the unique factorization property of $\intz[\x]$
that prevented us from finding a triple of polynomials in $\intz[\x]$
parametrizing \Pt s. Before we construct a parametrization
of \Pt s by a triple of integer-valued polynomials, we remark in
passing that $\Int(\intz^n)$ does not enjoy unique factorization
into irreducibles. An example of non-unique factorization into 
irreducibles in $\Int(\intz)$ is given by 
$$x(x-1)\ldots(x-k+1)=k! {x\choose k}.$$
The lefthand side is a product of $k$ irreducibles, while the 
righthand side, after factorization of $k!$ in $\intz$, becomes a
product of far more irreducibles (for large $k$).
(See \cite{CaCh97ivp} for integer-valued polynomials in general,
and \cite{AnCaChSm95FP,CahCha95El,ChaMcC05IPFE} for factorization
properties.)

\profess{Theorem.}
There exist $f,g,h\in \Int(\intz^4)$ such that $(f,g,h)$
parametrizes the set of \Pt s (as the variables range through $\intz$)
namely,
$$\displaylines{
\Biggl( {(2x-xw)((y+zw)^2-(z-yw)^2)\over 2},\hfill\cr
\hfill \vphantom{\Biggl)} (2x-xw)(y+zw)(z-yw), \hfill\cr
\hfill \vphantom{\Biggl(} {(2x-xw)((y+zw)^2+(z-yw)^2)\over 2}\Biggr).\cr
}$$
\endprofess

\proof
Every \Pt\ $(x,y,z)$ with $\gcd(x,y,z)=1$ and $z>0$ is either of the
form
$$T_1(a,b)=(a^2-b^2,\; 2ab,\; a^2+b^2),$$
or of the form
$$T_2(a,b)=(2ab,\; a^2-b^2,\; a^2+b^2),$$
with $a,b\in \intz$.
Since $$ 2\, T_2(a,b) = T_1(a+b,\> a-b),$$ every \Pt\ with
$\gcd(x,y,z)=1$ and $z>0$ is of
the form $c\, T_1(a,b)/2$ with $c\in\{1,2\}$ and $a,b\in \intz$.
Let 
$$T(a,b,c)=
\left({c(a^2-b^2)\over 2},\; cab,\; {c(a^2+b^2)\over 2}\right).$$
Then every \Pt\ is of the form $T(a,b,c)$ with 
$a,b,c\in\intz$. Also, every triple $T(a,b,c)$ with $a,b,c\in\intz$
is a rational solution of $x^2+y^2=z^2$.

So, the set of \Pt s is precisely the set of integer triples in
the range of the function $T\colon \intz^3\rightarrow \ratq^3$.

Now $T(a,b,c)\in \intz^3$ if and only if $c\congr 0$ mod $2$ or
$a\congr b$ mod $2$.
Triples $(a,b,c)\in\intz^3$ satisfying this condition can be 
parametrized by (for instance) 
$$(y+zw,\; z-yw,\; 2x-xw).$$ 
Indeed, if $w$ is even then $c\congr 0$ mod $2$, if $w$ is odd then 
$a\congr b$ mod $2$, and all $(a,b,c)$ satisfying either congruence
actually occur for some $(x,y,z,w)\in\intz^4$, as can be seen by
setting $w=0$ or $w=1$.

Therefore, substituting $y+zw$ for $a$, $z-yw$ for $b$,
and $2x-xw$ for $c$ in $T(a,b,c)$ yields a parametrization of
the set of \Pt s by a triple of integer-valued polynomials.
\endproof

\profess{Remark.}
The set of positive \Pt s is parametrized by
$$\displaylines{
\Biggl( {(x+(1-w)^2x)((y+(1+w)z)^2-y^2)\over 2} \vphantom{\Biggr)},\hfill\cr
\hfill\vphantom{\Biggl({1\over f}} (x+(1-w)^2x)(y+(1+w)z)y, \hfill\cr
\hfill \vphantom{\Biggl(} {(x+(1-w)^2x)((y+(1+w)z)^2+y^2)\over 2}\Biggr).\cr
}$$
where $x,y,z$ range through the positive integers and $w$ through the
non-negative integers. From this formula, a parametrization of positive
\Pt s with integer parameters can be obtained (using the 4-square theorem)
by replacing $w$ by $w_1^2+w_2^2+w_3^2+w_4^2$ and $x,\,y,\,z$ by
$x_1^2+x_2^2+x_3^2+x_4^2+1$, $y_1^2+y_2^2+y_3^2+y_4^2+1$, 
and $z_1^2+z_2^2+z_3^2+z_4^2+1$, respectively.
\endprofess

\proof
As in the proof of the theorem above, the positive \Pt s are
precisely the triples with positive integer coordinates in the
range of the function $T\colon \intz^3\rightarrow \ratq^3$.
Now $T(a,b,c)$ is a positive triple if and only if $a,b,c$ are
positive integers with $a> b$ and either $c\congr 0$ mod $2$ or
$a\congr b$ mod $2$.
Such triples $(a,b,c)$ are parametrized by (for instance)
$$(y+(1+w)z,\> y,\> x+(1-w)^2x)$$
with $x,y,z>0$ and $w\ge 0$. Therefore substituting $y+(1+w)z$ for
$a$, $y$ for $b$ and $x+(1-w)^2x$ for $c$ in $T(a,b,c)$ gives a
parametrization of positive \Pt s where $w$ ranges through 
non-negative integers and $x,y,z$ through positive integers.
The 4-square theorem allows us to convert this to a parametrization
with 16 integer parameters.
\endproof


\bigskip
\line{\bf References\hfil}%
\bibliography{pyth}
\bibliographystyle{siamese}

\bigskip\bigskip
\goodbreak
\line{
\vbox{
\hbox{S.F.\hfil}
\hbox{Institut f\"ur Mathematik C\hfil}
\hbox{Technische Universit\"at Graz\hfil}
\hbox{Steyrergasse 30\hfil}
\hbox{A-8010 Graz, Austria\hfil}
\tt
\hbox{frisch@blah.math.tu-graz.ac.at\hfil}
}
\hskip1.5cm
\vbox{
\hbox{L.V.\hfil}
\hbox{Department of mathematics\hfil}
\hbox{Pennsylvania State University\hfil}
\hbox{University Park, PA 16802\hfil}
\hbox{U.S.A.\hfil}
\tt
\hbox{vstein@math.psu.edu\hfil}
}
\hfill

}

\bye